\documentclass[doublespace,10pt]{amsart}
\usepackage{amsmath,amssymb,amsfonts,amscd}
\theoremstyle{plain}

\newtheorem{fthm}{Th\'eor\`eme}
\newtheorem{fprop}[fthm]{Proposition}

\newtheorem{fcor}[fthm]{Corollaire}

\newtheorem{prop}{Proposition}

\newtheorem{thm}[prop]{Theorem}
\newtheorem{coro}[prop]{Corollary}

\newcommand{\Pic}{{\rm Pic}}
\newcommand{\Hom}{{\rm Hom}}

\newcommand{\SL}{{\rm SL}}
\newcommand{\GL}{{\rm GL}}

\newcommand{\sC}{{\mathcal C}}

\newcommand{\sL}{{\mathcal L}}

\newcommand{\sO}{{\mathcal O}}

\newcommand{\C}{{\Bbb C}}

\newcommand{\G}{{\Bbb G}}
\renewcommand{\H}{{\Bbb H}}

\newcommand{\Z}{{\Bbb Z}}

\newcommand{\an}{{\rm an}}
\begin{document}
\title{Monodromies of algebraic connections on the trivial bundle}
\author{Byungheup Jun}
\email{byungheup.jun@uni-essen.de}
\address{FB 6, Mathematik, 
Universit\"at Essen,
45 117 Essen,
Germany}
\date{25. June, 2000}
\maketitle

\begin{abstract}
In this note, we study monodromies of algebraic connections
on the trivial vector bundle. We prove that on a smooth complex
affine curve, any monodromy arises as the underlying local system of an
algebraic connection on the trivial bundle. We give a
generalization of this for rank one monodromies in higher dimension.
\\[9pt]
{\bf\normalsize Monodromies des connexions alg\'ebriques sur le fibr\'e 
trivial}\\[5pt] 
{\sc R\'esum\'e.} Dans cette note, nous \'etudions la monodromie
de connexions alg\'ebriques sur le fibr\'e vectoriel trivial.
Nous prouvons que sur une courbe affine, lisse, complexe, toute
monodromie est sous-jacente \`a une connexion alg\'ebrique sur
le fibr\'e trivial. Nous donnons une g\'en\'eralisation de ceci
aux monodromies de rang un en dimension sup\'erieure.

\end{abstract}

\section*{Version fran\c caise abr\'eg\'ee}

Soit $U$ une vari\'et\'e affine lisse sur $\C$. Deligne a
montr\'e dans \cite{D} qu'il y a une correspondence biunivoque
entre les syt\`emes locaux sur $U$ et les fibr\'es vectoriels
alg\'ebriques \`a connexion plate, r\'eguli\`ere singuli\`ere
\`a l'infini. Si on supprime la condition de r\'egularit\'e \`a
l'infini, alors il y a beaucoup de connexions plates alg\'ebriques
\`a monodromie fix\'ee. Dans \cite{DEL}, il prouve que si $\dim U=1$, 
toute monodromie de rang 1 
est sous-jacente \`a une connexion alg\'ebrique  sur
le fibr\'e trivial. D\`es lors, les singularit\'es \`a l'infini
peuvent \^etre irr\'eguli\`eres.

\begin{fthm}[Deligne \cite{DEL}]
Sur une courbe affine lisse complexe $U$,  tout caract\`ere $\rho$ 
du groupe fondamental $\pi_1(U_\an)$ est sous-jacent \`a une
connexion alg\'ebrique sur le fibr\'e trivial.

\end{fthm}

{\it Id\'ee de la preuve.}--
On identifie 
$H^1_{dR}(U)=H^1(U_\an,\C)$ (\cite{Gd}), ce qui associe \`a la classe
de de~Rham de $\omega$ dans $\Gamma(U,\Omega^1_U)$ son
exponentielle  
$\exp[\omega]$ dans $H^1(U_\an,\C^*)$, 
qui correspond \`a la monodromie sous-jacente \`a la connexion
$(\sO,d-\omega)$.
L'annulation de  $H^2(U_\an,\Z)$ permet de conclure.
\hfill{\bf c.q.f.d.}

Nous donnons une version l\'eg\`rement diff\'erente du
th\'eor\`eme 1 dans le corollaire suivant.

\begin{fcor}
Soit $U$ une vari\'et\'e affine lisse complexe.
Les conditions suivantes sont \'equivalentes.
\begin{itemize}
\item[1)] $H^2(U_\an,\Z)$ n'a pas de torsion.
\item[2)] 
Tout caract\`ere $\rho$ 
du groupe fondamental $\pi_1(U_\an)$ est sous-jacent \`a une
connexion alg\'ebrique sur le fibr\'e trivial.
\item[3)] Soit $\sL$ un fibr\'e alg\'ebrique de rang 1 avec classe
de Chern triviale dans $H^2_{dR}(U)$. Alors
tout caract\`ere $\rho$ 
du groupe fondamental $\pi_1(U_\an)$ est sous-jacent \`a une
connexion alg\'ebrique sur $\sL$.
\end{itemize}  
\end{fcor}

Dans le th\'eor\`eme suivant, nous donnons un exemple
d'application. 
\begin{fthm}
Soit $X$ une vari\'et\'e projective lisse complexe et soit $D$
un diviseur ample supportant le groupe de
N\'eron-Severi  $NS(X)$ de $X$.
Alors tout caract\`ere $\rho$ 
du groupe fondamental $\pi_1(U_\an)$ est sous-jacent \`a une
connexion alg\'ebrique sur le fibr\'e trivial.
\end{fthm}

Nous renvoyons au diagramme (5) de Remark 1 dans \S{2} pour une
visualisation des r\'esultats pr\'ec\'edents.

Soit $U$ de nouveau une courbe lisse affine complexe.
Nous g\'en\'eralisons le th\'eor\`eme 1
aux repr\'esentations de rang sup\'erieur. Nous partons
de la d\'ecomposition d'une repr\'esentation $\rho$ \`a valeurs dans
 $\GL(r)$ comme produit tensoriel
$\chi\otimes\phi$, o\`u $\chi$ est un caract\`ere 
et
$\phi$ est une repr\'esentation \`a valeurs dans $SL(r)$.

Nous consid\'erons la suite exacte
\begin{equation}
0 \to \mu_r \to  \G_m\times\SL(r) \stackrel{\pi}{\to} \GL(r) \to 0,
\end{equation}
o\`u 
$\mu_r$ est le groupe des racines $r$-i\`emes de  l' unit\'e, 
et o\`u la premi\`ere application est d\'efinie par 
$\zeta \mapsto (\zeta, (1/\zeta)I)$ et la seconde est simplement
la multiplication d'une  matrice par un scalaire.

Alors on a la suite exacte longue de cohomologie sur $U_\an$
(d'ensembles point\'es)
\begin{equation}
\cdots \to H^1(U_\an,\G_m\times\SL(r)) \stackrel{\pi_*}{\to}
H^1(U_\an,\GL(r)) \stackrel{\delta}{\to} H^2(U_\an,\mu_r) \to \cdots
\end{equation}

Le groupe $H^2(U_\an,\mu_r)$ est nul car on est en dimension 1.
Donc  $\pi_*$ est
surjective. Par cons\'equent, il suffit de prouver le
th\'eor\`eme pour des repr\'esentations \`a d\'eterminant
trivial. 
\begin{fprop}
Toute repr\'esentation  $\rho: \pi_1(U_\an)\to SL(r,\C)$ 
du groupe fondamental est sous-jacente \`a une
connexion alg\'ebrique sur le fibr\'e trivial.
\end{fprop}
{\it Id\'ee de la preuve.--}
Soit $\rho: \pi_1(U_\an)\to SL(r,\C)$  et soit $(E,\nabla)$ l'unique
fibr\'e \`a connexion r\'eguli\`ere singuli\`ere \`a l'infini.
Sur une courbe affine $U$,  un fibr\'e alg\'ebrique $E$ 
a une d\'ecomposition 
$\sO^{r-1}\oplus\sL$ (cf. \cite{At}, \cite{Bou}).
Le d\'eterminant $\det(E)= \sL$ a la connexion d\'eterminant,
qui est aussi  
r\'eguli\`ere singuli\`ere \`a l' infini. Comme
${\rm det \ }\rho = 1$, on a $({\rm det}( E), {\rm det 
}(\nabla))=(\sO_U,d)$, donc en particulier
$\det(E)=\sL=\sO$. 
\hfill{\bf c.q.f.d.}

Finalement, nous obtenons
\begin{fthm}
Soit $U$ une courbe affine lisse complexe. Alors,
toute repr\'esentation  $\rho: \pi_1(U_\an)\to GL(r,\C)$ 
du groupe fondamental est sous-jacente \`a une
connexion alg\'ebrique sur le fibr\'e trivial.
\end{fthm}

\section{Introduction}
On a nonsingular projective algebraic variety $X$ over $\C$,
using the existence theorem of Cauchy-Kovalevski and Serre's 
GAGA(\cite{GAGA}), 
we have the 
Riemann-Hilbert correspondence: there exists a unique algebraic 
vector bundle with a 
connection for a given representation of the fundamental group $\pi_1(X)$
and vice versa.

On a noncomplete variety, 
the uniqueness no longer holds: there exist several vector bundles with 
connection for a monodromy.
In \cite{D}, Deligne showed the existence and uniqueness with regularity 
condition at ininity.

Instead of regularity, in this article, we think about a restriction 
of the underlying vector bundle of the connection.
In \cite{DEL}, 
Deligne showed one can realize any monodromy 
representation of the fundamental group of an affine curve $U$ 
as the underlying monodromy of 
a connection $\nabla=d-\omega$ on the trivial bundle of rank
1 on $U$(ie. $\sO_U$).

The main theorem(Theorem 5) generalizes Deligne's result(Theorem 1):
On an affine curve, to any monodromy representation of the fundamental 
group, there is a connection on trivial bundle with the prescribed underlying 
monodromy.

The aim of this note is to introduce Deligne's result on rank 1 case 
in \cite{DEL} and to generalize it in two directions: one is when $U$ is of 
higher dimension and the other is the similar to the higher rank 
representation of the fundamental group of the affine curve. 

{\bf Acknowledgements.}
Theorem 1, on which this note relies, is essentially due to Pierre Deligne.(\cite{DEL}).
I thank to H\'el\`ene Esnault for all the fruitful discussions on of my 
work and for sharing her time and ideas.

\section{Rank 1 representations}

\begin{thm}[Deligne \cite{DEL}]
On a smooth affine curve $U$, 
to given rank 1 monodromy representation $\rho$ of $\pi_1(U_\an)$,
there exists an algebraic connection $\nabla$ on $\sO$
with the underlying monodromy $\rho$.
\end{thm}

{\bf Proof.} 
From the exponential sequence, we have the following
exact sequence of singular cohomology groups.
\begin{equation}
H^1(U_\an,\C) \stackrel{\exp}{\longrightarrow} H^1(U_\an,\C^*) \longrightarrow 
H^2(U_\an,\Z) \longrightarrow H^1(U_\an, \C),
\end{equation}
where $U_\an$ denotes $U$ with the strong topology.
We know $H^1(U_\an,\C^*)$ is the group of isomorphism classes 
of rank 1 local systems
on $U_\an$
and $H^2(U_\an,\Z)$ vanishes by dimension reason. 

Via the identification(cf. \cite{Gd})
$$
H^*_{dR}(U):= \H^*(U,\Omega^\bullet) = H^1(U_\an,\C),
$$
the de~Rham class of $\omega$ in $H^0(U,\Omega^1)$ 
yields a class $[\omega]$ in $H^1(U_\an,\C)$, such that 
$$
\exp[\omega] 
\in H^1(U_\an,\C^*)=\Hom(\pi_1(U_\an),\C^*)
$$
is the underlying monodromy of $(\sO, d-\omega)$.

On an affine variety, the higher coherent sheaf cohomology vanishes.
Thus $H^0(U,\Omega^1)$ generates $H^1_{dR}(U)$.
Therefore $H^0(U,\Omega^1)$ surjects onto $\Hom(\pi_1(U_\an),\C^*)$.
This finishes the proof.
\hfill {\bf  q.e.d.}\\[5pt]

Without much difficulty, we generalize the theorem when $\dim U > 1$.

\begin{coro}
Let $U$ be an affine variety over $\C$. 
Then the followings are equivalent:
\begin{itemize}
\item[1)] $H^2(U_\an,\Z)$ is torsion-free.
\item[2)] To any rank 1 representation $\rho$ of $\pi_1(U_\an)$, there exists an algebraic 
 connection $\nabla$ on $\sO$ with underlying monodromy $\rho$.
\item[3)] If a line bundle $\sL$ has vanishing Chern class in $H^2_{dR}(U)$.
To any rank 1 representation $\rho$ of $\pi_1(U_\an)$ and , there exists an algebraic connection
$\nabla$ on $\sL$ with underlying monodromy $\rho$.
\end{itemize}  
\end{coro}
{\bf Proof. }
1) $\Leftrightarrow$ 2)
As in the proof of the theorem 1, one has
the exact sequence arising from the exponential sequence(3).
On an affine variety $U$, $\Gamma(U,\Omega^1_{\rm closed})$ generates 
$H^1_{dR}(U) = H^1(U_\an,\C)$. 
Then $\exp$ in the exact sequence is surjective if and only if 
$\Gamma(U,\Omega^1_{\rm closed})\to H^1(U_\an,\C^*)$ is surjective or equivalently,
$H^2(U_\an,\Z)$ is torsion-free.

2) $\Rightarrow$ 3)
Consider the complex of sheaves on the Zariski topology of $U$:
$$
\sC^\bullet : 
	\sO^* \stackrel{d\log}{\longrightarrow} \Omega^1 \to \Omega^2 \to \cdots.
$$
The first hypercohomology of $\sC^\bullet$ is the group of isomorphism classes of line 
bundles with integrable connection, which will be denoted by 
$\Pic^\nabla(U)$.
The filtration of $\sC$ by its degree gives
\begin{equation}
\begin{CD}
\H^1(U,\sC^\bullet) @>>> H^1(U,\sO^*)	@>>> 
		\H^2(U,\Omega^{\bullet \ge 1})	\\
 		@|				@|		@|    \\
\Pic^\nabla(U) @>>>  \Pic(U)		@>c>> H^2_{dR}(U),
\end{CD}
\end{equation}
where $c$ denotes the first Chern class map and the last vertical map is 
an isomorphism since
the kernel and the cokernel are $H^1(U,\sO)=H^2(U,\sO)=0$.

Since $c(\sL)=0$, from the exact sequence, there exists an
integrable connection $\nabla_0$ on $\sL$.
Let $\rho_0$ be the monodromy of $(\sL,\nabla_0)$.
As one can realize any monodromy on $\sO$, one has an algebraic connection
$\nabla_1$ with the monodromy $\rho_0^{-1}\cdot\rho$ on $\sO$.
Twisting $(\sL, \nabla_0)$ with $(\sO,\nabla_1)$, one has the connection
$\nabla=\nabla_0\otimes\nabla_1$ on $\sL$ with the desired monodromy 
$\rho$.

3) $\Rightarrow$ 2) Clear.
\hfill{\bf q.e.d.}

As a special case of the corollary, we have the following:
\begin{thm}
Let $X$ be a smooth projective variety over $\C$ and $D$ be an ample divisor supporting
the N\'eron-Severi group $NS(X)$ of $X$.
Then to any monodromy $\rho$ in $\Hom(\pi(U_\an),\C^*)$, there exists an integrable algebraic 
connection $\nabla=d-\omega$ on $\sO_U$ with underlying monodromy $\rho$.   
\end{thm}

{\bf Proof.}
By the equivalence of the corollary 2, we have to show that $H^2(U_\an,\Z)$ is 
torsion-free.
But, in the exact sequence (3), torsion comes from $H^1(U_\an,\C^*)$.
Using the existence theorem in \cite{D}, 
one has a line bundle on $X$ with connection $(\sL,\nabla)$ 
which is regular singular at $D$ with the underlying monodromy $\rho$. 
The image of $\sL$ maps into $NS(X)$, and thus a torsion element in $H^2(U_\an,\Z)$
vanishes. 
Therefore $H^2(U_\an,\Z)$ is torsion-free.
\hfill{\bf q.e.d.}\\[5pt]
{\bf Remark 1.}
Using a commutative diagram, one can summarize the result in a simple form 
where all rows and columns are exact:

\begin{equation}
\begin{array}{ccccccccc}
  &					& 0 & & 0 & & 0 \\
  &					& \downarrow   &	 & \downarrow & & \downarrow \\
0 & \longrightarrow & \displaystyle{\frac{H^0(U,\Omega^1)_\Z}{d \log H^0(U,\sO^*)}} & \longrightarrow 
  					& \displaystyle{\frac{H^0(U,\Omega_{\rm closed}^1)}{d \log H^0(U,\sO^*)}} & \longrightarrow & H^1(U_\an,\C^*) & \longrightarrow & 0\\
  &	   	 	   	 	& \downarrow     & 				  &	\downarrow &  & \downarrow \\
0 & \longrightarrow & \Pic^{\nabla_0}(U) & \longrightarrow & \Pic^\nabla(U) & 
  					\longrightarrow & \Hom(\pi_1(U_\an),\C^*) & \longrightarrow & 0\\
  & 	& \downarrow & 	   	   &  \downarrow &	  &  \downarrow \\
0 & \longrightarrow	& \Pic(U)	&	\longrightarrow	 &	    \Pic(U)	& \longrightarrow & 0	 \\
  &		& \downarrow &	 		& \downarrow \\
  &  & H^2_{dR}(U) 	&	=		&	H^2_{dR}(U) &    
\end{array}
\end{equation}
In the diagram, 
$\Pic^{\nabla_0}(U)$ is the subgroup of $\Pic^\nabla(U)$ with respect to
the trivial monodromy and $H^0(U,\Omega^1)_\Z$ is the global 1-forms
with integral periods.
When $\dim U=1$, one has $H^2_{dR}(U)=0$. 
Thus, we have $\Pic^\nabla(U)$ and $\Pic^{\nabla_0}(U)$ as extensions 
of the $\Pic(U)$.

\section{Higher rank case on curves}
In this section, $U$ is again a smooth affine curve over $\C$.

One has the short exact sequence:
\begin{equation}
\begin{split}
0 \to \mu_r  \to \G_m\times\SL(r) & \stackrel{\pi}{\to} \GL(r) \to 0,\\
        \zeta \mapsto   (\zeta, (1/\zeta) I) &\\
                         (\ell,M)  &\mapsto \ell M
\end{split}
\end{equation}
where $\mu_r$ denotes the multiplicative group of $r$-th roots of
unity.
Because $\mu_r$ maps into the center of $\G_m \times \SL(r)$,
we have long exact sequence of cohomologies on $U_\an$ (of pointed sets).
\begin{equation}
\cdots \to H^1(U_\an,\G_m\times\SL(r)) \stackrel{\pi_*}{\to} H^1(U_\an,\GL(r)) 
\stackrel{\delta}{\to}
H^2(U_\an,\mu_r) \to \cdots
\end{equation}
By dimension reason, $H^2(U_\an,\mu_r)$ vanishes and one sees $\pi_*$ 
surjective.
Thus it suffices to solve the problem for $SL$-representations.

If $A$ is a Dedekind domain and $M$ is a projective $A$-module of finite 
rank $r$,
one has $M\cong A^{r-1}\oplus L$ where $r= rank(M)$ and $L \in \Pic(A)$.
As a special case, a vector bundle $E$ of $rank(E)=r$ on an affine curve $U$
is isomorphic to $\sO^{r-1}\oplus \sL$, where $\sL = \det E$ is a line bundle 
on $U$ (cf. \cite{Bou}, \cite{At})

\begin{prop}
Let $U$ be a smooth affine curve.
To any $\SL(r)$-representation $\rho$ of $\pi_1(U_\an)$, 
there exists an algebraic connection $\nabla$
on the trivial bundle of rank $r$ with underlying monodromy $\rho$.
\end{prop}
{\bf Proof.} 
Let $\rho$ be a $\SL(r)$-representation and $X=U\cup D$ be the smooth completion of $U$. 
Then there exists a vector bundle $E$ on $X=U\cup D$ the smooth completion 
of $U$, with a connection $\nabla$ which is regular singular, 
$$
\nabla : E \to E\otimes \Omega^1(\log D),
$$
such that the restriction on $U$ gives the monodromy $\rho$ (cf. \cite{D}).

By the decomposition theorem, one has $E|_U = \sO_U^{r-1}\oplus \sL$, 
where $\sL = \det(E|_U)$.
On $X$, 
$$
\det \nabla : \det(E) \to \det(E)\otimes \Omega^1(\log D)
$$
is regular at $D$ and has trivial monodromy since the representation  
is $SL$-valued.
Using the uniqueness theorem in \cite{D},
we have 
$$
(\det (E)|_U,\det \nabla) = (\sO_U,d).
$$ 
In particular $\sL=\sO_U$ and $E_U = \sO^r_U$.
\hfill {\bf q.e.d.}

After the proposition, we obtain the main theorem.

\begin{thm}
Let $U$ be a smooth affine curve over $\C$.
Then to any rank $r$ representation of 
the fundamental group $\rho\in \Hom(\pi_1(U_\an),\GL(r))$,
there exists a connection on the trivial bundle of rank $r$ on an affine 
open set $U$ with underlying monodromy $\rho$.
\end{thm}
{\bf Proof.} Clear.
\hfill{\bf q.e.d.}\\[5pt]
\bibliographystyle{plain}
\renewcommand\refname{References}

\end{document}